\documentclass[12pt]{amsart}
\usepackage{amssymb,amsmath,amsthm}
\usepackage{enumerate}

\DeclareMathOperator{\Ran}{Ran}

\DeclareMathOperator{\Ker}{Ker}
\DeclareMathOperator{\Tr}{Tr}

\DeclareMathOperator*{\slim}{s-lim}
\DeclareMathOperator{\rank}{rank}
\DeclareMathOperator{\iindex}{index}
\DeclareMathOperator{\modd}{mod}
\DeclareMathOperator{\flow}{sp.flow}

\newcommand{\abs}[1]{\lvert#1\rvert}

\newcommand{\norm}[1]{\lVert#1\rVert}

\newcommand{\R}{{\mathbb R}}

\newcommand{\calH}{{\mathcal H}}

\newcommand{\calF}{\mathcal{F}}
\newcommand{\calN}{\mathcal{N}}

\numberwithin{equation}{section}

\theoremstyle{plain}
\newtheorem{theorem}{\bf Theorem}[section]

\theoremstyle{definition}

\theoremstyle{remark}
\newtheorem*{remark*}{\bf Remark}

\newcommand{\ess}{\text{\rm ess}}
\newcommand{\ac}{\text{(ac)}}

\begin{document}

\title[Integer-valued Birman-Krein formula]
{An integer-valued version of the Birman-Krein formula}

\author{Alexander Pushnitski}
\address{Department of Mathematics,
King's College London, 
Strand, London, WC2R~2LS, U.K.}
\email{alexander.pushnitski@kcl.ac.uk}

\dedicatory{To the memory of M.~Sh.~Birman}

\begin{abstract}
We discuss an identity in abstract scattering theory which 
can be interpreted as an integer-valued version of the 
Birman-Krein formula.  
\end{abstract}



\maketitle

\section{Introduction}\label{sec.a}

In the classic paper \cite{BK} by M.~Sh.~Birman and M.~G.~Krein, 
the following identity (known today as the Birman-Krein formula)
was proven:
\begin{equation}
e^{-2\pi i\xi(\lambda)}=\det S(\lambda). 
\label{a1}
\end{equation}
Here $S(\lambda)$ and $\xi(\lambda)$ are the scattering matrix and 
the spectral shift function of a pair of self-adjoint operators
in a Hilbert space  
and $\lambda\in\R$ is a spectral parameter;
precise definitions will be given in Sections~\ref{sec.b} and \ref{sec.c}.
The identity \eqref{a1}
first appeared in \cite{Lifshits2,BuslaevFaddeev,Buslaev66}
in some important particular cases; 
in \cite{BK} it was recognised as an abstract theorem 
of mathematical scattering theory. 

The Birman-Krein formula is invariant with respect to adding 
an integer to $\xi(\lambda)$. Thus, one can say that 
\eqref{a1} involves only the fractional part of $\xi(\lambda)$; 
the information about the integer part of $\xi(\lambda)$ is 
``lost''. 
The purpose of this note is to state and discuss an identity 
(see \eqref{e4} below)
which can be interpreted as an integer-valued version 
of the Birman-Krein formula. 

In Sections~\ref{sec.b} and \ref{sec.c}
we give the necessary definitions.  
The main result is stated in Section~\ref{sec.e}. 
Proofs will appear in \cite{Push5}.

The author is grateful to D.~Yafaev, T.~Suslina, B.~Davies and N.~Filonov 
for useful discussions
and remarks on the text of the paper.

\section{$\xi(\lambda)$ and $\Xi(\lambda)$}\label{sec.b}

\subsection{Naive definition}
Throughout the paper, $A$ and $B$ are self-adjoint operators
in a Hilbert space $\calH$ such that the difference 
$$
V=B-A
$$
is a compact operator. 
For simplicity of exposition we will also assume that $A$ and $B$
are bounded.  
The assumption of compactness of $V$ implies, in particular, 
that the essential spectra of $A$ and $B$ coincide: 
$\sigma_\ess(A)=\sigma_\ess(B)$. 
We denote by $E_A(\lambda)$ the spectral projection of $A$
corresponding to the interval $(-\infty,\lambda)$ 
and let 
$$
N_A(\lambda):=\rank E_A(\lambda)=\Tr E_A(\lambda)\leq \infty
$$
be the eigenvalue counting function of $A$. 
For $\lambda<\inf\sigma_\ess(A)$ the difference
\begin{equation}
N_A(\lambda)-N_B(\lambda)
\label{b1}
\end{equation}
is well defined, since both $N_A(\lambda)$ and $N_B(\lambda)$ are finite. 
The difference \eqref{b1} measures the shifts of the eigenvalues
of $B$ relatively to the eigenvalues of $A$. 
For $\lambda>\inf\sigma_\ess(A)$ the difference \eqref{b1} 
formally gives $\infty-\infty$. 
Below we discuss two regularisations of \eqref{b1}:
the spectral shift function $\xi(\lambda)$ and the index function $\Xi(\lambda)$. 

\subsection{The spectral shift function $\xi(\lambda)$}

Assume that the difference $V=B-A$ is a trace class operator. 
Then \cite{Lifshits1,Krein} the following Lifshits-Krein trace formula 
holds true: 
\begin{equation}
\Tr(\varphi(B)-\varphi(A))
=
\int_{-\infty}^\infty \varphi'(t)\xi(t)dt, 
\quad \forall \varphi\in C_0^\infty(\R). 
\label{b2}
\end{equation}
Here $\xi(\cdot)=\xi(\cdot;B,A)$ is a uniquely defined
function in $L^1(\R)$ which is called the 
spectral shift function (SSF). 
See \cite{BYa} or 
\cite[Section~8]{Yafaev} for a detailed exposition 
of the SSF theory. 

Formally taking $\varphi=\chi_{(-\infty,\lambda)}$ 
($=$ the characteristic function of $(-\infty,\lambda)$), 
we obtain 
\begin{equation}
\xi(\lambda;B,A)
=
\Tr(E_A(\lambda)-E_B(\lambda))
=
\Tr E_A(\lambda) - \Tr E_B(\lambda)
=
N_A(\lambda)-N_B(\lambda), 
\label{b3}
\end{equation}
whenever the r.h.s. makes sense. 
In particular, this calculation is not difficult to justify for $\lambda<\inf\sigma_\ess(A)$.
It shows that $\xi(\lambda)$ 
is the natural regularisation of the difference 
$N_A(\lambda)-N_B(\lambda)$.

\subsection{The index of a pair of projections}\label{sec.b3}
Let $P,Q$ be orthogonal projections in a Hilbert space;
consider the spectrum of the difference $P-Q$. 
It is obvious that $\sigma(P-Q)\subset[-1,1]$. 
By using the commutation relation 
$$
W(P-Q)=-(P-Q)W, \quad W=I-P-Q,
$$
it is not difficult to see 
\cite[Theorem~4.2]{ASS}
that 
\begin{equation}
\dim \Ker(P-Q-\lambda I)=\dim \Ker(P-Q+\lambda I), 
\quad \lambda\not=\pm 1.
\label{b5}
\end{equation}
A pair $P,Q$  is called Fredholm, if
\begin{equation}
\{1,-1\}\cap\sigma_{\rm ess}(P-Q)=\varnothing.
\label{b6}
\end{equation}
In particular, if $P-Q$ is compact, then the pair $P$, $Q$ is Fredholm.
The index of a Fredholm pair is defined by the formula
\begin{equation}
\iindex(P,Q)=\dim\Ker(P-Q-I)-\dim\Ker(P-Q+I).
\label{b7}
\end{equation}
We note that $\iindex(P,Q)$ coincides with the Fredholm index
of the operator $QP$ viewed as a map from $\Ran P$ to $\Ran Q$, 
see \cite[Proposition~3.1]{ASS}.

If $P-Q$ is a trace class operator,  then 
\begin{equation}
\iindex(P,Q)=\Tr(P-Q),
\label{b8}
\end{equation}
since all the eigenvalues of $P-Q$ apart from $1$ and $-1$  in the series 
$\Tr(P-Q)=\sum_k\lambda_k(P-Q)$
cancel out by \eqref{b5}. 
It follows that in the simplest case of finite rank projections $P,Q$ we have
\begin{equation}
\iindex(P,Q)=\rank P-\rank Q.
\label{b9}
\end{equation}

\subsection{The index function $\Xi(\lambda)$}\label{sec.b4}
Suppose that for some $\lambda\in\R$, the pair $E_A(\lambda)$, $E_B(\lambda)$ 
is Fredholm. 
Then we will say that the index $\Xi(\lambda)=\Xi(\lambda;B,A)$ 
exists and define it by 
\begin{equation}
\Xi(\lambda;B,A)
=
\iindex(E_A(\lambda),E_B(\lambda)).
\label{b10}
\end{equation}
The function $\Xi(\lambda;B,A)$ 
has already appeared in the literature in various guises 
(see e.g. \cite{APS,ADH,DH,Hempel2,Sobolev,RobbinSalamon,GMN,GM,
BPR,ACS,KMS,Hempel,Simon100DM});
its properties  were reviewed  and proven in a systematic fashion 
in \cite{Push3}.

For $\lambda<\inf\sigma_\ess(A)$, both projections 
$E_A(\lambda)$, $E_B(\lambda)$ have finite rank and so 
by \eqref{b9} we have
$$
\Xi(\lambda;B,A)
=
N_A(\lambda)-N_B(\lambda), 
\quad
\lambda<\inf\sigma_\ess(A).
$$
Thus, $\Xi(\lambda)$, along with $\xi(\lambda)$,  
is a natural regularisation of the difference 
$N_A(\lambda)-N_B(\lambda)$. 
Using the Riesz integral  representation for the 
spectral projections and the resolvent identity,
it is not difficult to prove   that 
for all $\lambda\in\R\setminus\sigma_\ess(A)$, the difference
$E_A(\lambda)-E_B(\lambda)$ is compact and therefore 
$\Xi(\lambda)$
exists.
Below we will also give a criterion for the existence
of $\Xi(\lambda)$ on the essential spectrum of $A$ and  $B$, 
see Theorem~\ref{thm0}.

\subsection{Comparison of  $\xi$ and $\Xi$}\label{sec.b5}
\begin{enumerate}[1.]
\item
If $V$ is a trace class operator, we have 
$$
\Xi(\lambda)=\xi(\lambda),
\quad 
\forall \lambda\in\R\setminus\sigma_\ess(A). 
$$
In particular, this holds true for all $\lambda\in\R$ if $\dim\calH<\infty$. 
\item
$\Xi(\lambda)$ is always integer-valued, whereas $\xi(\lambda)$ 
is, in general, real-valued when $\lambda$ belongs to the essential spectrum of $A$, $B$.  
\item
If $V$ is trace class then the SSF $\xi(\lambda)$ is automatically 
well defined for a.e. $\lambda\in\R$. 
The function $\Xi(\lambda)$ may not exist on a set of positive measure
even if $V$ is a rank one operator; see \cite{KM}.
On the other hand, the existence of $\Xi(\lambda)$ does not require
trace class assumptions.  
\end{enumerate}

\section{The scattering matrix}\label{sec.c}

\subsection{Strong smoothness assumptions}\label{sec.c1}
Below we make assumptions typical for smooth scattering theory,
which goes back to \cite{Faddeev} and \cite{Kato66}.
We fix a compact interval $\Delta=[a,b]\subset\R$ and assume that 
the spectrum of $A$ in $\Delta$ is purely absolutely continuous
with  constant multiplicity $N\leq\infty$. 
In the terminology of \cite{Yafaev}, we
assume that the operator $G=\abs{V}^{1/2}$ is strongly $A$-smooth 
on $\Delta$ with some exponent $\gamma\in (0, 1]$.
This means the following. 
Let $\calF$ be a unitary operator from $\Ran E_A(\Delta)$ to 
$L^2(\Delta,\calN)$, $\dim\calN=N$, such that $\calF$
diagonalizes $A$: if $f\in\Ran E_A(\Delta)$ then 
\begin{equation}
(\calF  A f)(\lambda)
=
\lambda (\calF f)(\lambda), 
\quad \lambda\in\Delta.
\label{a8}
\end{equation}
The strong $A$-smoothness of  $G$ 
on the interval $\Delta$ means that the operator 
\[
G_\Delta\overset{\rm def}{=}GE_A(\Delta) : \Ran E_0(\Delta) \to \calH
\] 
satisfies the condition
\begin{equation}
(\calF G_\Delta^*\psi)(\lambda)
=
Z(\lambda)\psi, 
\quad 
\forall \psi\in\calH,
\quad
\lambda\in\Delta,
\label{a3b}
\end{equation}
where $Z=Z(\lambda):\calH\to\calN$
is a family of compact operators obeying
\begin{equation}
\norm{Z(\lambda)}\leq C,
\quad
\norm{Z(\lambda)-Z(\lambda')}\leq C\abs{\lambda-\lambda'}^\gamma,
\quad
\lambda, \lambda'\in\Delta.
\label{a3a}
\end{equation}

\subsection{The scattering matrix}\label{sec.c2}
Under the above strong smoothness assumption, 
the local wave operators
$$
W_\pm
=
W_\pm(A,B;\Delta)
=
\slim_{t\to\pm\infty}e^{itB}e^{-itA}E_A(\Delta)
$$
exist and are complete, 
i.e. $\Ran W_+=\Ran W_-=\Ran E_B^{\ac}(\Delta)$; 
here $E_B^{\ac}(\cdot)$ is the absolutely continuous 
part of the spectral measure of the operator $B$. 
The local scattering operator ${\mathbf S}=W_+^*W_-$ is unitary 
in $\Ran E_A(\Delta)$ and commutes with $A$. 
Thus, we have a representation
\[
( \calF \mathbf{S}  \calF^* f)(\lambda) =S(\lambda) f(\lambda), 
\quad \text{ a.e. }\lambda\in \Delta,
\]
where the operator  $S(\lambda):\calN\to\calN$  
is  called the scattering matrix
for the pair of operators $A$, $B$. 
The scattering matrix  is a unitary operator in $\calN$.
Under the above assumptions, one can prove that 
$S(\lambda)-I$ is a compact operator for all 
$\lambda\in\Delta$.
Thus, the spectrum 
of $S(\lambda)$ consists of eigenvalues on the unit circle 
$e^{i\theta_n(\lambda)}$ with the only possible point of accumulation 
being  $1$. 
Further, one can prove that $S(\lambda)$ depends
continuously on $\lambda$ in the operator norm. 
For the details, see the original papers \cite{Faddeev,Kato66} or
the survey \cite{BYa2} or the book \cite{Yafaev}.

\subsection{The spectral flow of $S(\lambda)$}\label{sec.c3}
Let us recall the definition of  
the spectral flow of the family  $\{S(\lambda)\}_{\lambda\in[a,b]}$.  
The spectral flow is an integer-valued function $\mu$ 
on $\mathbb T\setminus\{1\}$. The naive definition of the spectral flow is 
\begin{multline}
\flow(e^{i\theta}; \{S(\lambda)\}_{\lambda\in[a,b]})=
\\
\langle \text{the number of eigenvalues of $S(\lambda)$ which cross
$e^{i\theta}$ in the anti-clockwise direction}\rangle
\\
-
\langle \text{the number of eigenvalues of $S(\lambda)$ which cross
$e^{i\theta}$ in the clockwise direction}\rangle,
\label{d4}
\end{multline}
as $\lambda$ increases monotonically from $a$ to $b$. 
Here $\theta\in(0,2\pi)$ and the eigenvalues
are counted with multiplicities taken into account. 
The eigenvalues of $S(\lambda)$ may cross $e^{i\theta}$ infinitely many times, 
and thus the above naive definition needs to be replaced by a more rigorous 
one. Below we describe one such possible definition; there are other approaches
to this definition in the literature, see e.g. \cite{APS,RobbinSalamon}. 
 
Let us introduce some notation for the eigenvalue counting function of 
$S(\lambda)$. 
For $\theta_1,\theta_2\in(0,2\pi)$ we denote 
$$
N(e^{i\theta_1},e^{i\theta_2}; S(\lambda))
=
\sum_{\theta\in[\theta_1,\theta_2)}
\dim\Ker(S(\lambda)-e^{i\theta}I), 
\quad \text{ if }
\theta_1<\theta_2,
$$
and
$$
N(e^{i\theta_1},e^{i\theta_2}; S(\lambda))
=
-N(e^{i\theta_2},e^{i\theta_1}; S(\lambda))
\quad \text{ if }
\theta_1>\theta_2. 
$$
Assume first that there exists $\theta_0\in(0,2\pi)$ 
such that $e^{i\theta_0}\notin\sigma(S(\lambda))$ for all $\lambda\in[a,b]$. 
Then one can define the spectral flow of the family 
$\{S(\lambda)\}_{\lambda\in[a,b]}$  by 
\begin{equation}
\flow(e^{i\theta};\{S(\lambda)\}_{\lambda\in[a,b]})
=
N(e^{i\theta},e^{i\theta_0};S(b))
-
N(e^{i\theta},e^{i\theta_0};S(a)).
\label{d5}
\end{equation}
It is evident that this definition is independent of the choice of $\theta_0$ and 
agrees with the naive definition \eqref{d4} whenever the latter makes sense.

In general, $\theta_0$ as above may not exist. 
However, by a standard argument based on the compactness of $[a,b]$ 
one can always find 
the values $a=\lambda_0<\lambda_1<\dots<\lambda_n=b$ 
such that 
for each of the subintervals $\Delta_i=[\lambda_{i-1},\lambda_i]$, 
a point $\theta_0$  with the required properties
can be found. 
Thus, the spectral flow of each of the corresponding families
$\{S(\lambda)\}_{\lambda\in\Delta_i}$ is well defined. Now one can set
\begin{equation}
\flow(e^{i\theta};\{S(\lambda)\}_{\lambda\in[a,b]})
=
\sum_{i=1}^n
\flow(e^{i\theta}; \{S(\lambda)\}_{\lambda\in\Delta_i}).
\label{d5a}
\end{equation}
It is not difficult to see that the above definition is independent 
of the choice of the subintervals $\Delta_i$ and agrees with the naive 
definition \eqref{d4}.

\section{Main results}\label{sec.e}

\subsection{Statement of the results}\label{sec.e1}
The first preliminary result concerns the existence of $\Xi(\lambda)$ on the continuous
spectrum. 
\begin{theorem}\cite{Push0}\label{thm0}
Assume that for some interval $\Delta=[a,b]$ 
the spectrum of $A$ in $\Delta$ is purely absolutely continuous and let
the strong smoothness assumption \eqref{a3b}, \eqref{a3a} hold true. 
Then for all $\lambda\in(a,b)$ one has
$$
\sigma_{\ess}(E_B(\lambda)-E_A(\lambda))
=
[-\alpha(\lambda),\alpha(\lambda)], 
\quad
\alpha(\lambda)=\tfrac12\norm{S(\lambda)-I}. 
$$
In particular, $\Xi(\lambda;B,A)$ exists if and only if 
$-1\notin\sigma(S(\lambda))$. 
\end{theorem}
A description of the absolutely continuous spectrum 
of the difference $E_B(\lambda)-E_A(\lambda)$
is also available in terms of the spectrum of $S(\lambda)$, 
see \cite{PushYa}.

As $\lambda$ increases monotonically in the interval $\Delta$, 
the eigenvalues of $S(\lambda)$ rotate on the unit circle and the quantity 
$\alpha(\lambda)$ changes continuously in $\lambda$. 
According to Theorem~\ref{thm0}, 
the index $\Xi(\lambda)$ exists if and only if $\alpha(\lambda)<1$, 
i.e. if and only if the spectrum of $S(\lambda)$ does not contain 
the point $-1$.

Further simple analysis based on the stability of Fredholm index 
shows that the function $\Xi(\lambda)$ is constant on the intervals 
where $-1\notin\sigma(S(\lambda))$. 
Thus, 
the  integer-valued function $\Xi(\lambda)$ can jump only 
at the points $\lambda$ where $-1\in\sigma(S(\lambda))$.
This leads to the natural question: what is the size of the 
jump of $\Xi(\lambda)$ when an eigenvalue of $S(\lambda)$ 
crosses $-1$? 
Our main result below answers this question. 

\begin{theorem}\label{th1}
Assume that for some interval $\Delta=[a,b]$ 
the spectrum of $A$ in $\Delta$ is purely absolutely continuous and let
the strong smoothness assumption \eqref{a3b}, \eqref{a3a} hold true. 
Fix $\lambda_1,\lambda_2\in(a,b)$, $\lambda_1<\lambda_2$ and 
assume that $-1\notin\sigma(S(\lambda_1))$ and $-1\notin\sigma(S(\lambda_2))$. 
Then 
\begin{equation}
\Xi(\lambda_2; B,A)
-
\Xi(\lambda_1; B,A)
=
-\flow(-1; \{S(\lambda)\}_{\lambda\in[\lambda_1,\lambda_2]}).
\label{e1}
\end{equation}
\end{theorem}
The proof will appear in \cite{Push5}. 
The proof is based, roughly speaking, on 
a continuous  deformation of the pair of operators
$A$, $B$, into a pair of operators which has 
an ``infinitesimal spectral gap'' at a point $\lambda\in(a,b)$.
This deformation, together with the Birman-Schwinger
principle, enables us  to calculate 
$\Xi(\lambda;A,B)$ in terms of some auxiliary operators.  
This calculation makes it possible to relate 
$\Xi(\lambda;A,B)$ to the 
eigenvalue counting function $N(-1,e^{i\theta_0};S(\lambda))$. 
From this relation it is not difficult to derive \eqref{e1}.

Theorem~\ref{th1} can be extended to a certain 
class of unbounded operators $A$, $B$. 
In many concrete examples of interest, the 
absolutely continuous spectra 
of $A$ and $B$ contain the semi-axis $[0,\infty)$, 
and $\norm{S(\lambda)-I}\to0$ as $\lambda\to\infty$. 
In this case one can take $\lambda_2\to\infty$ in \eqref{e1}, which yields
\begin{equation}
\Xi(\lambda;B,A)
=
\flow(-1; \{S(\lambda')\}_{\lambda'\in[\lambda,\infty)}). 
\label{e4}
\end{equation}
This is discussed in \cite{Push5} for the case 
$A=-\Delta$, $B=-\Delta+V$ in $\calH=L^2(\R^d)$, 
 with a short range potential $V$.

\subsection{Comparison of \eqref{e4} with the Birman-Krein formula}\label{sec.e2}
Let $\{e^{i\theta_n(\lambda)}\}$ be the eigenvalues of $S(\lambda)$ distinct from $1$. 
To compare \eqref{e4} with the Birman-Krein formula \eqref{a1}, let 
us rewrite the latter as 
\begin{equation}
\xi(\lambda)
=
-\frac1{2\pi}
\sum_n \theta_n(\lambda)
\quad (\modd \  1).
\label{e3}
\end{equation}
As discussed in Section~\ref{sec.b}, 
the left hand sides of \eqref{e4} and \eqref{e3} 
are two different regularisations of $N_A(\lambda)-N_B(\lambda)$. 
The right hand sides of \eqref{e4} and \eqref{e3} are two different
quantities related to the spectrum of the scattering matrix.

The identity \eqref{e4} relates two integers. 
The identity \eqref{e3} relates two real numbers modulo 1. 
Thus, in some sense (perhaps yet to be understood) they 
present complementary pieces of information.

\end{document}